\documentclass{amsart}[12pt]
\usepackage{amssymb}
\usepackage{mathrsfs}
\usepackage{epsfig}

\input xy
\xyoption{all}

\newcommand{\A}{\ensuremath{\mathbb A}}

\newcommand{\q}{\ensuremath{\mathfrak q}}

\newcommand{\m}{\ensuremath{\mathfrak m}}

\newcommand{\NN}{\ensuremath{\mathbb N}}
\newcommand{\PP}{\ensuremath{\mathbb P}}
\newcommand{\QQ}{\ensuremath{\mathbb Q}}
\newcommand{\RR}{\ensuremath{\mathbb R}}
\newcommand{\ZZ}{\ensuremath{\mathbb Z}}
\newcommand{\CC}{\ensuremath{\mathbb C}}
\newcommand{\OO}{\ensuremath{\mathcal O}}
\newcommand{\II}{\ensuremath{\mathcal I}}

\newcommand{\la}{\ensuremath{\lambda}}
\newcommand{\w}{\ensuremath{v}}
\newcommand{\ga}{\ensuremath{\mathfrak a}}
\newcommand{\vl}{\ensuremath{{\rm vol}}}

\newcommand{\hh}[3]{\ensuremath{h^{#1}\left(#2,#3\right)}}
\newcommand{\HH}[3]{\ensuremath{H^{#1}\left(#2,#3\right)}}
\newcommand{\hff}{\ensuremath{\left( d\lfloor j(i)l \rfloor +j(i)\right)\overline{H} -\lfloor j(i)l \rfloor F}}
\newcommand{\si}[2]{\ensuremath{ \left( { \left( #1 \right) }^{#2}\right) } }
\newcommand{\egy}{\ensuremath{ I_1(k)} }
\newcommand{\ket}{\ensuremath{ I_2(k)} }

\newcommand{\fl}[1]{\ensuremath{ \lfloor #1 \rfloor }}
\newcommand{\zj}[1]{\ensuremath{ \left( #1 \right) }}

\theoremstyle{definition}

\newtheorem*{ex}{Example}

\theoremstyle{plain}

\newtheorem{thm}{Theorem}[section]
\newtheorem{prop}{Proposition}[section]
\newtheorem{cor}{Corollary}[section]

\begin{document}

\title{A divisorial valuation with irrational volume}
\author{Alex K\"uronya}

\maketitle

\section{Introduction}

The purpose of this paper is to construct a divisorial valuation with irrational volume. Let $(R,\m,{\CC})$ be an $n$-dimensional noetherian local ring and consider a rank one valuation $\w$ of its fraction field  centered on  $R$ (ie. $\w$ is nonnegative on $R$ and strictly positive on  $\m$). Then one can associate to $\w$ its volume
\begin{equation} \vl (\w)=\limsup_{m}{\frac{{\rm length}(R/{\q}_m)}{m^n/n!}}  \label{1} \ . \end{equation}
Here ${\q}_k$ denotes the ideal of elements with valuation at least $k$. 
This is an analogue of the Samuel multiplicity  
\begin{equation} e(\ga ):=\limsup_{m}{\frac{{\rm length}(R/{\ga}^m)}{m^n/n!}} \label{2} \end{equation}
of an $\m$-primary ideal $\ga\subseteq R$. In fact, if ${\q}_m={\ga}^m$ for a fixed ideal $\ga$ then it is evident that $e(\ga )=\vl (\w)$. The volume of a valuation has implicitly been studied already in \cite{cs}, but it was first explicitly defined in \cite{els}. The terminology is intended to emphasize the relation with global invariants of linear series on projective varieties. 

A natural question is to what extent the properties of $\vl(\w)$ mirror those of the Samuel multiplicity. A result of \cite{els} and \cite{mus} asserts that 
\begin{equation} \vl(\w) =\lim_{m\rightarrow\infty}{\tfrac{e({\q}_m)}{m^n}} \label{3} \end{equation}
so in any event the volume is governed by multiplicity. A basic fact about multiplicity --- which is not apparent from the definition above --- is that $e(\ga)$ is always an integer. However, simple examples show that this is false for the volume:

\begin{ex}
Consider the monomial valuation of ${\CC}(x,y)$ centered at the origin of ${\A}^2$ and  defined by $\w (x)=1, \w(y)=\alpha$ for $\alpha\in \RR$. It follows directly from the definition that 
\begin{equation} \vl(\w)= \frac{1}{\alpha} \label{4} \end{equation}  which is irrational if $\alpha$ is.
\end{ex}

In the example however,  the irrationality was `built into' the valuation in the sense that if one chooses $\alpha$ to be rational then the volume will also be rational. On the other hand, suppose $\w$ is a divisorial valuation, that is, a  valuation with $ rk\ {\w}=1$ and   $tr.deg_{\CC}\ \w=n-1$ (or, in other words, a valuation attached to an irreducible exceptional divisor of a birational map). In dimension two, Cutkosky and Srinivas (\cite{cs}, Corollary 1) prove that under mild hypotheses the volume is indeed rational. Their proof relies on the existence of Zariski decompositions on surfaces. Our objective here is to show that in higher dimensions there are divisorial valuations with irrational volume. 

Specifically, we prove
\begin{thm} Let $R={\CC}[x_1,x_2,x_3,x_4]_{(x_1,x_2,x_3,x_4)}$. There exists a  divisorial valuation $\w$ of  $\CC (x_1,x_2,x_3,x_4)$ centered in $R$ such that \[ \vl(\w)\not\in \QQ .\]
\end{thm}
\noindent
Note that divisorial valuations always have value group $\ZZ$. A related invariant of a valuation is the associated graded ring \[ {gr}_{\w}R:={\bigoplus}_{m\geq 0}\frac{{\q}_m}{{\q}_m^+} \] where ${\q}_m^+=\{ g\in K|\w (g)>m\}$. It is easily observed (for a proof see for example \cite{els}) that if ${gr}_{\w}R$ is finitely generated then $\vl(\w)$ is rational. As a consequence we obtain a simple construction of a divisorial valuation whose associated graded ring is not finitely generated (another example of this phenomenon has been described in \cite{cgp}, Proposition 2).

The  construction uses in particular ideas suggested by \cite{cs} and recent  irrationality results on asymptotic invariants of algebraic varieties (\cite{cut},\cite{cel}). One starts with a smooth curve $C\subseteq {\PP}^{3}$ with irrational asymptotic Castelnuovo--Mumford regularity. Realizing ${\PP}^{3}$ as the exceptional divisor of ${\mathrm Bl}_0({\CC}^4)$ the order of vanishing along $C$ determines a divisorial valuation on ${\CC}^4$. We relate the asymptotic regularity of $C$ to the volume of the resulting valuation to arrive at the desired conclusion.

{\bf Acknowledgments.} I would like to thank my advisor, Rob Lazarsfeld, for drawing my attention to this problem, for his support and many useful comments. I am also grateful to Mihnea Popa, Mike Roth, Jessica Sidman and Karen Smith for helpful discussions.

\section{A volume formula for certain divisorial valuations}

In this section we  give an explicit formula for the volume of divisorial valuations of a certain kind. We will consider  valuations $\w$ of the field $\CC (x_1,\dots ,x_n)$  centered at the origin $o$ of ${\A}^n$, hence the local ring will be $R={\OO}_{ {\A}^n,o}$. 

The divisorial valuations in question will be constructed by two successive blowups as follows: we start with the blowup $\pi :W=Bl_o({\A}^n)\rightarrow {\A}^n$ of the origin with exceptional divisor $V\simeq {\PP}^{n-1}$. Next, we pick a smooth subvariety $T\subseteq V$ and   form the blow-up $p:Y=Bl_T(W)\rightarrow W$ of $W$ along $T$ with exceptional divisor $E$. We denote the composition $p\circ \pi$ by $f$. The valuation $\w$ is the valuation determined by $E$; hence its valuation ring is ${\OO}_{Y,E}$. 

\medskip

\xymatrix{ & & & E  \ar[d]  & \ar@{} [r]-{\subseteq}   & Y=Bl_T(W) \ar[d]^-{p}  \ar@/^5pc/  [dd]_f \\
	   & & & T \ar@{} [r] |-{\subseteq}   & {\PP}^{n-1}=V  \ar[d]   \ar@{} [r] |-{\subseteq} & W=Bl_o({\A}^n) \ar[d]^-{\pi} \\
 & & &  & o \ar@{} [r] |-{\in} &  {\A}^n 
}
 
\medskip

\noindent
As ${\q}_m=f_*{\mathcal O}_Y(-mE)\subseteq {\mathcal O}_{{{\A}^n},o}$ is ${\m}_o$-primary, $ {\rm length}(R/{\q}_m)=dim_{\CC}{{\OO}_{{\A}^n,o}/{\q}_m}$ is finite. 

We will obtain an explicit formula for the colengths of the valuation ideals ${\q}_m$ in terms of the cohomology of  the ideal sheaf ${\II}_T\subseteq {\OO}_{{\PP}^{n-1}}$ of $T$ in ${\PP}^{n-1}$.

\begin{prop} With notation as above,
\begin{equation} {\rm length}(R/{\q}_m)=\sum_{s=0}^{m-1}{\left( h^0({\PP}^{n-1},{\OO}_{{\PP}^{n-1}}(s))-h^0({\PP}^{n-1},{\II}^{m-s}_T(s)) \right)}. \label{5} \end{equation}
\end{prop}

\noindent
{\bf Proof:\ } First observe that $h\in f_*{\OO}_Y(-mE)\subseteq {\CC}[x_1,\dots ,x_n]$ if and only if  ${\pi}^*h$ vanishes on $T$ to order at least $m$. Since the order of vanishing on $T$ is a local invariant, we can make computations in local coordinates. Specifically, in suitable local coordinates on an affine open subset $U\subseteq W$ meeting $T$, $\pi$ is given by

\xymatrix{ & & & & & (x_1,x_2,\dots ,x_n)\in U \ar[d]^-{{\pi}|_U} \ar@{} [r] |-{\subseteq} & W \ar[d]^{\pi}  \\
	   & & & & & (x_1,x_1x_2,\dots x_1x_n) \ar@{} [r] |-{\in}   & {\A}^n 
.} 
\noindent
Hence  if $h=\sum_{i_1,i_2,\dots ,i_n}{a_{i_1i_2\dots i_n}x_1^{i_1}\dots x_n^{i_n}}$ then

\begin{eqnarray} \left({\pi}^*h\right)(x_1,x_2,\dots ,x_n) & = & h(x_1,x_1x_2,\dots ,x_1x_n) \nonumber \\ 
 & = & \sum_{d}{{x_1}^d\sum_{i_2,\dots ,i_n}{a_{d-i_2-\dots -i_n,i_2,\dots ,i_n}{x_2}^{i_2}\dots {x_n}^{i_n}}} \\
 & = &  \sum_d{{x_1}^dg_d(x_2,\dots ,x_n)}\ .\nonumber
\label{7}
\end{eqnarray}
Here $g_d\in \CC [x_2,\dots ,x_n]$ is a degree $d$ polynomial and for $d_1\not=d_2$ the set of coefficients $a_{i_1i_2\dots i_n}$ involved in the polynomials $g_{d_1}$ and $g_{d_2}$ are disjoint. Therefore the conditions we get on the vanishing of various derivatives of the $g_{d_i}$'s are independent for different $d_i$'s.
We compute the partial derivatives of ${\pi}^*h$ on $T\subseteq \left\{x_1=0\right\}=V|_U$:

\begin{equation} \left( {\partial}^{m_1}_{x_1}{\partial}^{m_2}_{x_2}\dots {\partial}^{m_n}_{x_n} \right) \zj{ {\pi}^*h }(0,x_2,\dots ,x_n) = m_1!\left({\partial}^{m_2}_{x_2}\dots {\partial}^{m_n}_{x_n}g_{m_1} \right)(x_2,\dots ,x_n)\ .
\label{8} \end{equation}
For  ${\pi}^*h$ to vanish on $T$ up to order $m$ is the same as asking  
${\partial}^{m_2}_{x_2}\dots {\partial}^{m_n}_{x_n}g_{m_1}$ to vanish  identically  on $T$ for all $m_2+\cdots m_n<m-m_1$. This happens exactly if each $g_s$ vanishes to order $m-s$ on $T$, therefore  $g_s$ determines an element in \HH{0}{{\PP}^{n-1}}{{\II}_T^{m-s}(s)}. Hence we can deduce that the codimension of ${\q}_m=f_*{\OO}_Y(-mE)$ is 

\begin{equation} \sum_{s=0}^{m-1}{\left( h^0({\PP}^{n-1},{\OO}_{{\PP}^{n-1}}(s))-h^0({\PP}^{n-1},{\II}^{m-s}_T(s)) \right)}. \label{9} \end{equation}

\qed

\begin{cor} Let $\w$ be a divisorial valuation of $\CC (x_1,\dots ,x_n)$ centered at the origin of ${\A}^n$ as in the construction above. Then 

\begin{equation} \begin{array}{rl}
\vl (\w) &=\limsup_m {\tfrac{1}{m^n/n!}\sum_{s=0}^{m-1}{\left( {h^0({\PP}^{n-1},{\OO}_{{\PP}^{n-1}}(s))-h^0({\PP}^{n-1},{\II}^{m-s}_T(s))} \right)}} \\ \\
&=1-\liminf_m \tfrac{1}{m^n/n!}\sum_{s=0}^{m-1}{\hh{0}{{\PP}^{n-1}}{{\II}^{m-s}_T(s)}}
. \end{array} \label{10} \end{equation}
\end{cor}

\section{An example of a divisorial valuation with irrational volume}

Based on the  construction  of the previous section we will exhibit an example of a divisorial valuation with irrational volume. As said earlier, this example  also establishes  the existence of a divisorial valuation whose  associated graded ring  is not finitely generated. The source of irrationality is the choice of a certain configuration $C\subseteq S\subseteq {\PP}^{n-1}$ with $C$ having irrational asymptotic regularity (for the basic results on asymptotic Castelnuovo--Mumford regularity the reader can consult \cite{pag}, Section 1.8). The known examples of asymptotic irregularity involve either K3 or abelian surfaces. The instance we will use is the K3 surface constructed in \cite{cut}.

Using  \cite{mor}, Theorem 2.9., Cutkosky shows that there exists a K3 surface $S$ such that ${\rm Pic}(X)\simeq {\ZZ}^3$ and its intersection form is  $q(x,y,z)=4x^2-4y^2-4z^2$ in suitable coordinates. We choose this surface  $S$ for our computations. Then $S\subseteq {\PP}^3$ is a degree four surface, whose nef cone and the effective cone are equal and given by 
\[ {\rm  Nef}(S)=\{\alpha\in NS(S)_{\RR}| \left( {\alpha}^2\right)\geq 0,(\alpha\cdot h)\geq 0\} \]
where $h$ is  any ample class on $S$. Even though $S\subseteq {\PP}^{3}$ has degree four, for the clarity of the exposition we will write $d$ for its degree throughout the paper.

Fix a very ample divisor $H$ on $S$ that embeds it into ${\PP}^{3}$ with degree $d$ and pick $C$ to be  an effective divisor such that the line $tH-C$ in the N\'eron--Severi space  intersects the boundary of the nef cone at the irrational value $\la$. Then  the asymptotic irregularity $\lambda$ of $C$ with respect to the fixed very ample divisor $H$ will be  irrational (\cite{cut}). In the computations we will choose $H=(1,0,0)$ and $C=(9,1,1)$ on our K3 surface.

 The main ingredients  of the volume formula of the previous section are  the dimensions of the cohomology groups $\HH{0}{{\PP}^{3}}{{\II}_C^r(m)}$ which we will relate  to the cohomology of certain divisors on the blowup of $C\subseteq {\PP}^{3}$.  Let ${H}^{\prime}\subseteq {\PP}^{3}$ be a hyperplane such that ${H}^{\prime}\cdot S=H$ and denote the exceptional divisor of the blowup $\pi :X=Bl_C{\PP}^{3}\rightarrow {\PP}^{3}$ by $F$. Then $F={\pi}^{-1}C$ and  ${\pi}^*S=\overline{S}+F$ with $\overline{S}$  the strict transform of $S$. In what follows let $\overline{H}={\pi}^*{H}^{\prime}$.

One has  ${\pi}_*{\OO}_X(m\overline{H}-rF)\simeq {\II}^r_C(m)$ and $R^i{\pi}_*{\OO}_X(m\overline{H}-rF)=0$ for $i>0$ (cf. \cite{mat}, Proposition 10.2), and hence
\begin{equation} \hh{i}{{\PP}^{3}}{{\II}_C^r(m)}=\hh{i}{X}{m\overline{H}-rF} \label{12} \end{equation}
for all $i,m,r\geq 0$. The dimensions of the cohomology groups appearing in (\ref{12}) can be computed explicitly for $i=0$ thanks to Riemann--Roch. We will then interpret the sum of the  dominant terms of the $\hh{0}{{\PP}^{3}}{{\II}_C^r(m)}$'s appearing in the volume formula as a Riemann sum for a certain integral.

\begin{prop} One has
\begin{equation} \vl (\w)=1-4\left( \int_{\overline{BA}}{\si{m\overline{H}-rF}{3}} + \int_{\overline{OB}}{\si{m\overline{H}-rF}{3}} \right)\ .\label{26} \end{equation}
The integrals are computed using the parameterizations
\begin{equation} \overline{BA}:\ {\gamma}_1(t)=(t,1-t),\ \tfrac{\la}{\la +1}\leq t\leq 1, 
\end{equation}
\begin{equation} \overline{OB}:\ {\gamma}_2(t)=\left( \tfrac{\la}{\la -d}((d+1)t-d),\tfrac{1}{\la -d}((d+1)t-d)\right),\ \tfrac{d}{d+1}\leq t\leq \tfrac{\la}{\la+1}
\end{equation}
 over the  piecewise linear curve $\overline{ABO}$  with $O$ the origin, $A=(0,1)$ and $B$ the intersection of the lines $m=\la r$ and $m+r=1$. \label{formula}
\end{prop}

\begin{center} {\sc Figure 1} \\  \label{fig1}
\epsfig{file=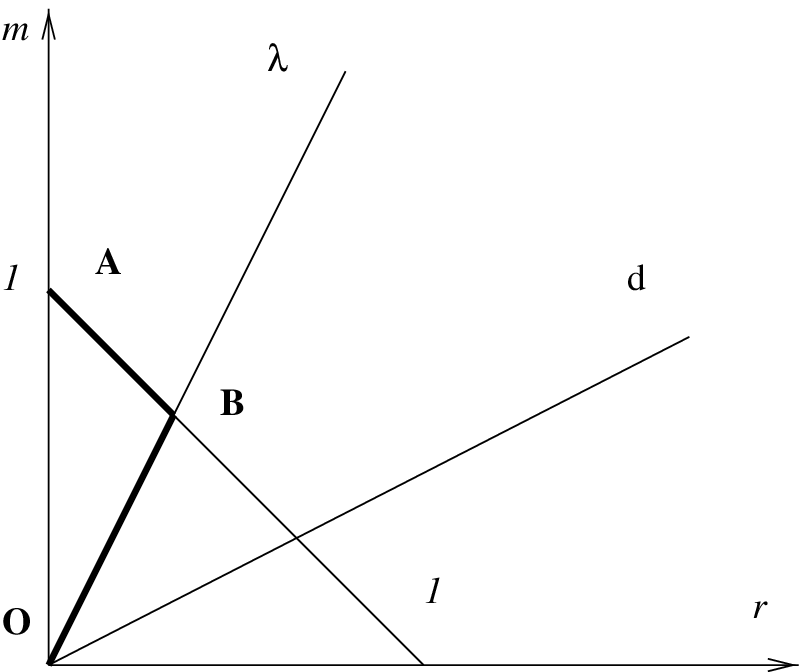,height=1.5in,width=2in}
\end{center}

The  integrals in the proposition are illustrated in Figure 1. We next give a detailed description of our results while postponing the proofs to the last section. First, we explain the computation of  $ \hh{0}{X}{m\overline{H}-rF}=\hh{0}{{\PP}^{3}}{{\II}^r_C(m)} $. There are three cases to the computation, depending on the ratio $\tfrac{m}{r}$. We will show that if $\tfrac{m}{r}>\la $ then
\begin{equation} \hh{0}{X}{m\overline{H}-rF}=	\chi (m\overline{H}-rF)\ .\label{19} 
\end{equation}
In the case $ \tfrac{m}{r} <\la$ we will prove 
\begin{equation} h^0(m\overline{H}-rF)=h^0((m-d)\overline{H}-(r-1)F)\ ,\label{22} \end{equation}
which used iteratively will either lead back to  the previous case --- if  $\tfrac{m}{r}>d$ --- or give  \[ h^0(m\overline{H}-rF)=0\] if $\tfrac{m}{r}<d$. This is illustrated on Figure 2 by the arrows between the dots. We summarize these results in the next proposition.

\begin{center} {\sc Figure 2} \\ 
\epsfig{file=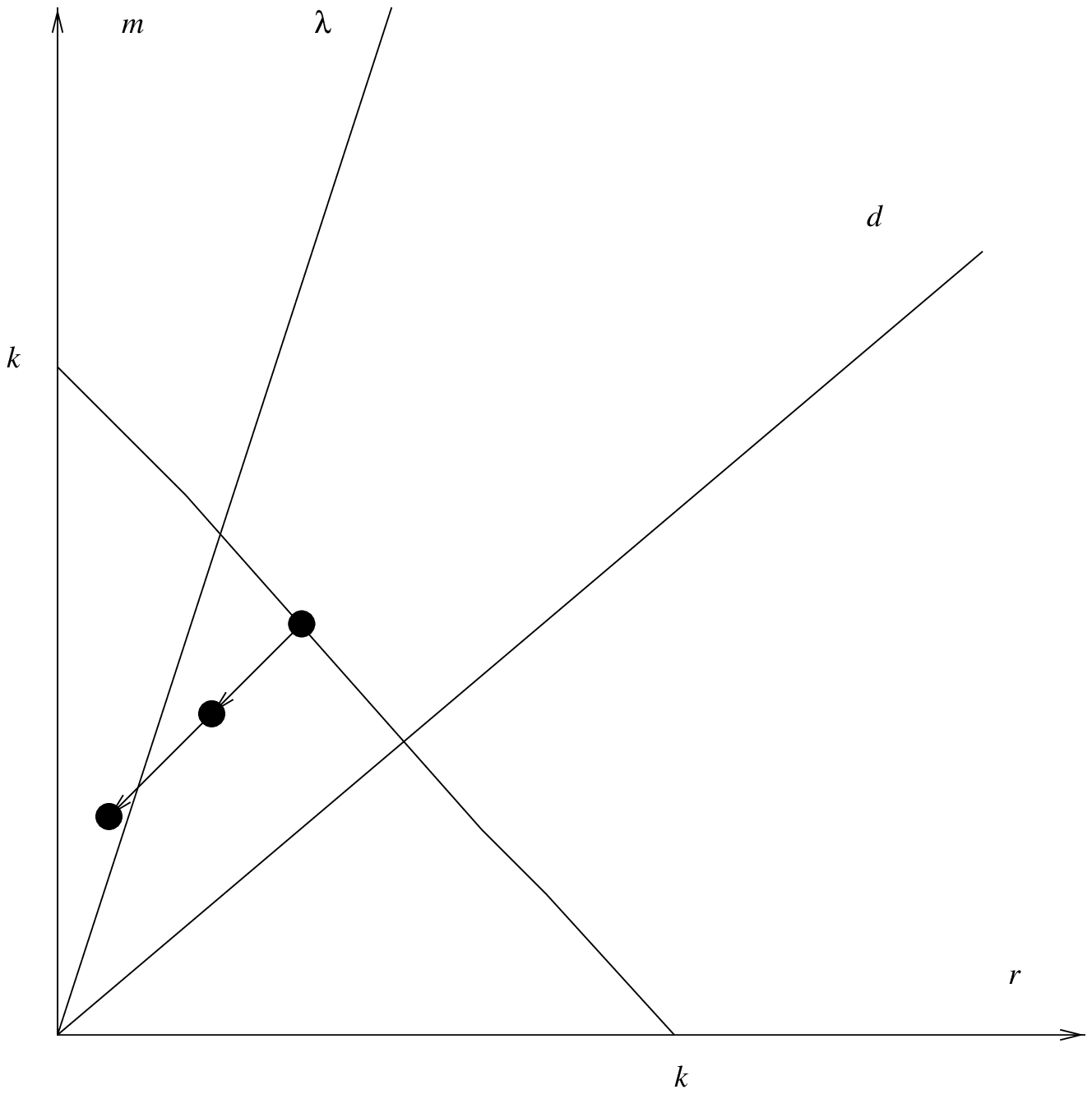,height=3in,width=3in} \label{fig2}
\end{center}

\begin{prop}  With notation as above,
\begin{equation}
\hh{0}{X}{m\overline{H}-rF}=
	\begin{cases}
	\chi (m\overline{H}-rF) & \text{if $\tfrac{m}{r}\geq \la$},\\
	\chi((d\lfloor\frac{m-dr}{\la -d}\rfloor+m-dr)\overline{H}-(\lfloor\frac{m-dr}{\la-d}\rfloor)F) & \text{ if $\la > \tfrac{m}{r}\geq d$},\\
	0  & \text{ if $d>\tfrac{m}{r}$.}
	\end{cases}
\label{18} 
\end{equation} \label{ezaz}
\end{prop}
\noindent
We obtain  the integral expression for the volume of $\w$  with these computations along with the volume formula of the previous section.

Using the formula of Proposition \ref{formula}  we can explicitly calculate the volume of the corresponding valuation (for example with the computer algebra package Maple).
The result 
\begin{equation} 
\vl (\w)=\frac{144629}{2352980}+\frac{1408}{588245}\sqrt{2}
\label{37} \end{equation}
\noindent
is indeed irrational.

\section{Proofs}

This section contains the proofs of Propositions 3.1 and 3.2. Before we move on to the proofs themselves, we make some observations. We keep the notation of the previous sections.

First,  by Kodaira vanishing on $S$ and the description of the effective cone one has
 \begin{equation} \hh{0}{S}{mH-rC}= \begin{cases}
\tfrac{1}{2}\chi\zj{S,mH-rC}  & \text{if\ } \la <\tfrac{m}{r} \\
0 & \text{otherwise.}
\end{cases}
\label{11} \end{equation}
Next, we set up a short exact sequence that we will use repeatedly. Specifically, tensoring the sequence
\begin{equation} 0\rightarrow {\OO}_X(-\overline{S})\rightarrow {\OO}_X\rightarrow {\OO}_{\overline{S}}\rightarrow 0 \label{13} \end{equation}
by ${\OO}_X((m+d)\overline{H}-(r+1)F)$ leads to 
\begin{equation}  0\rightarrow {\OO}_X(m\overline{H}-rF)\rightarrow {\OO}_X((m+d)\overline{H}-(r+1)F)\rightarrow  {\OO}_S((m+d)H-(r+1)C)\rightarrow 0  \label{egzakt} \end{equation}
for all $m,r\geq 0$ via the isomorphisms
\begin{equation} {\OO}_X(-d\overline{H}+F)\simeq {\OO}_X(-\overline{S})\label{15} \end{equation}
and
\begin{equation} {\OO}_S((m+d)H-(r+1)C)\simeq {\OO}_{\overline{S}}((m+d)\overline{H}-(r+1)F)\ .\label{16} \end{equation}

\noindent
We next prove Proposition \ref{ezaz}.\\
\noindent
{\bf Proof of \ref{ezaz}:}
According to \cite{cut}, Theorem 9, if $\tfrac{m}{r}>\la $ then all higher cohomology of $m\overline{H}-rF$ vanishes, so $\hh{0}{X}{m\overline{H}-rF}=\chi (m\overline{H}-rF)$. For  $\tfrac{m}{r}<\la $, consider the long exact sequence corresponding to  (\ref{egzakt}) with $m-d,r-1$ in the place of $m,r$:
\begin{equation} 0\rightarrow H^0(X,(m-d)\overline{H}-(r-1)F)\rightarrow H^0(X,m\overline{H}-rF)\rightarrow H^0(S, mH-rC)\rightarrow\cdots \ .\label{23} \end{equation}
Observe that if $\tfrac{m}{r}<\la$ then the last term is zero and  
\begin{equation}
 \hh{0}{X}{m\overline{H}-rF}=\hh{0}{X}{(m-d)\overline{H}-(r-1)F} \ .
\end{equation}
 We  can continue this process replacing again $m,r$ by $m-d,r-1$ until either $m-d<0$ which implies $\hh{0}{X}{(m-d)\overline{H}-(r-1)F}=0$ (this will happen exactly when $\tfrac{m}{r}<d$ for the starting pair) or  when $\tfrac{m-d}{r-1}$ becomes $\geq \la$. So far we have 
\begin{equation} \hh{0}{X}{m\overline{H}-rF}=\hh{0}{X}{m^{\prime}\overline{H}-r^{\prime} F} \label{24} \end{equation}
where $m^{\prime},r^{\prime}$ are of the form $m-ds$, $r-s$ ($s$ a  positive integer) such that either $m^{\prime}<0$ and hence $\hh{0}{X}{m\overline{H}-rF}=0$ or $m^{\prime}> r^{\prime}{\lambda}$. The former case happens if and only if $m<dr$.

In the latter case let $L$ be the line with slope $d$ and going through the point $(r,m)$. Then the integral point with the biggest $r$-coordinate in $L\cap \left\{ \tfrac{m}{r}>\la \right\} $ is $\left( r^{\prime},m^{\prime}\right)$. In concrete terms

\begin{equation} m^{\prime} =d\lfloor \frac{m-dr}{{\lambda}-d} \rfloor + (m-dr),\ r^{\prime} =\lfloor \frac{m-dr}{{\lambda}-d} \rfloor\ . \label{25} \end{equation}
As $\tfrac{m^{\prime}}{r^{\prime}} > \la $ (hence all  higher cohomology of $m^{\prime} \overline{H}-r^{\prime} F$ vanishes) this completes the proof.\\
\qed \\

\noindent
Finally, we move on to the proof of the integral formula for the volume.\\
\noindent
{\bf Proof of \ref{formula}:\ } According to the volume formula of Section 2.,
\begin{equation} \vl (\w) = 1-\varliminf_k \tfrac{4!}{k^4} \sum_{i=0}^{k-1}{\hh{0}{{\PP}^{3}}{{\II}^{k-i}_C(i)}}\ . \label{27}
\end{equation}
In terms of the $(r,m)$-plane, we add  the terms $\hh{0}{X}{m\overline{H}-rF}$ on the line segment $m+r=k\ (m,r\geq 0)$ for fixed $k$. First we introduce some notation. Let 
\begin{equation} j(i):=(d+1)i-dk \text{\ and\ } l:=\tfrac{1}{\la -d}\ ,\ \label{28} \end{equation}
\begin{equation} I_1(k):=\left\{ i |\ d<\tfrac{i}{k-i}<\la\right\} \ ,\label{29} \end{equation}
\begin{equation} I_2(k):=\left\{ i |\ \la \leq\tfrac{i}{k-i}\right\} \ .\label{30} \end{equation}
In other words, $I_1(k)$ is the part of the line segment $m+r=k\ (m,r\geq 0)$ that falls in  the region between the lines $\tfrac{m}{r}=d$ and $\tfrac{m}{r}=\la$ while $I_2(k)$ is the part  falling in the region between the $m$-axis and the line $\tfrac{m}{r}=\la$.

\begin{equation} \begin{array}{rl}
\sum_{i=0}^{k-1}{\hh{0}{{\PP}^{3}} {{\II}_C^{k-1}(i)}}&=  \sum_{i=0}^{k-1}{\hh{0}{X}{i\overline{H}-(k-i)F}}\\ \\
 &= \sum_{\egy}{ \chi\zj{ \hff}} \\ \\
 &+ \sum_{\ket}{ \chi\zj{ (i\overline{H}-(k-i)F)}}.
\end{array} \label{31}
\end{equation}
Riemann--Roch theorem on $3$-folds implies that 
\begin{equation} \chi (i\overline{H}-(k-i)F)=\frac{1}{3!}\left( (i\overline{H}-(k-i)F)^{3}) +O(k^{2}\right) \label{32} \end{equation}
and so

\begin{equation} \begin{array}{rl}
 \chi\zj{ \hff}  &= \frac{1}{3!}\si{\hff}{3} \\ \\ &+ O(k^{2})\ .
\end{array}\label{33} \end{equation}
Therefore

\begin{equation} \begin{array}{rl}
 \vl (\w) &= 1-\varliminf_k\bigg( \tfrac{4}{k^4} \sum_{\egy}{ \si{\hff}{3} } \bigg. \\
 &+ \bigg. \tfrac{4}{k^4} \sum_{\ket}{\si{ i\overline{H}-(k-i)F}{3}}\bigg)  \\ \\
 &= 1-\varliminf_k\bigg( \tfrac{4}{k} \sum_{\egy}{ \si{d\tfrac{\lfloor jl\rfloor +j}{k}\overline{H}-\tfrac{\lfloor jl\rfloor}{k}F}{3} } \bigg. \\
 &+ \bigg. \tfrac{4}{k} \sum_{\ket} {   \si{ \tfrac{i}{k}\overline{H}-\zj{1-\tfrac{i}{k}}F}{3}   }  \bigg)\ .
\end{array}
\label{34} \end{equation}
For $n\in\NN$ and  $x\in\RR$ one has $\vert x-\tfrac{1}{n}\fl{nx} \vert \leq \tfrac{1}{n}\label{35} $. As  $ \tfrac{j(i)}{k}=(d+1)\tfrac{i}{k}-d$ this then implies

\begin{equation}
 \begin{array}{rl}
 \vl (\w) &=  1-\varliminf_k\bigg( \tfrac{4}{k} \sum_{\egy}{ \si{\zj{dl+1}\zj{\zj{d+1}\tfrac{i}{k}-d}\overline{H}-l\zj{\zj{d+1}\tfrac{i}{k}-d}F}{3}} \bigg. \\
&+ \bigg. \tfrac{4}{k} \sum_{\ket}{\si{\tfrac{i}{k}\overline{H}-\zj{1-\tfrac{i}{k}}F}{3}}\bigg) \\ \\
&= 1 -4\left( \int_{\tfrac{d}{d+1}}^{\tfrac{\la}{\la +1}}{\si{\zj{\zj{dl+1}\zj{\zj{d+1}t-d}}\overline{H}-l\zj{\zj{d+1}t-d}F}{3} dt} \right. \\ \\
& \left. + \int_{\tfrac{\la}{\la +1}}^{1}{ \si{t\overline{H}-(1-t)F}{3} dt} \right) \ .
\end{array}
\label{36} \end{equation} \qed

\medskip\medskip

{\small\sc
\noindent
department of mathematics, university of michigan, ann arbor, mi 48109, usa \\
{\em email address:\ }{\tt akuronya@umich.edu} \\

\noindent
and \\

\noindent
computer and automation institute of the hungarian academy of sciences, h-1518 budapest, p.o. box 63., hungary

}

\end{document}